# An Application of Ptolemy's Theorem: Integral triangles with a 120° angle and with the bisector (of the 120° angle) also of integral length

Konstantine Zelator
Mathematics, Statistics, and Computer Science
212 Ben Franklin Hall
Bloomsburgh University of Pennsylvania
400 East 2$^{nd}$ Street
Bloomsburg, PA 17815
USA
and P.O. Box 4280
Pittsburgh, PA 15203
kzelator@bloomu.edu
e-mails: konstantine_zelator@yahoo.com

October 5, 2011

## 1    Introduction

During the academic year 2008-2009, while I was a mathematics faculty member at Rhode Island College, one of my colleagues, Professor Vivian LaFerla mentioned the following geometry problem to me (see Figure 1). Let *ABC* be an equilateral triangle and consider its circumscribed circle. Pick a point *P* on the arc $\overline{BC}$ which does not contain the vertex *A* and let *D* be the point of intersection between the line segments $\overline{AP}$ and $\overline{BC}$. Prove that $\dfrac{1}{PD} = \dfrac{1}{PB} + \dfrac{1}{PC}$ (here $\overline{PD}, \overline{PB}, \overline{PC}$ are the lengths of the line segments respectively). By inspection, we see that the angle ∠*BPC* has measure 120° and that *PD* is the angle bisector of the 120° angle. Accordingly, the above statement says that in a triangle with a 120° angle, the reciprocal of the length of the angle bisector of the 120° angle is the sum of the reciprocals of the lengths of the two sides which contain the 120° angle. We prove this in Proposition 1 in Section 3. To do so, we employ Ptolemy's Theorem for cyclic quadrilaterals, a very small known historical result in geometry. We state Ptolemy's Theorem in Section 2. References for Ptolemy's Theorem abound. For example, see [1], [2], and [3]. All the aforementioned is geometric in nature. However, the main objective of this work is number-theoretic in nature. Namely, to parametrically describe the set of all integral triangles with a 120° angle and whose angle bisectors (of their 120 degree angles) also have integral lengths. To achieve this objective, which we do in Section 7, two goals must be accomplished.

(1) The general solution, in positive integers, of the Diophantine equation
    $\dfrac{1}{z} = \dfrac{1}{x} + \dfrac{1}{y}$. We derive all the positive integer solutions to this equation
    in Propostition 2, Section 5. Toward that end we employ Lemma 1
    (Euclid's Lemma) which is stated in Section 4 (this is a very well-known lemma).

(2) The parametric description of the set of integral triangles with a 120° angle. An integral triangle is a triangle whose side lengths are integers. We state formulas that describe (parametrically) the set of all integral triangles with a 120° angle, in Section 6. This is work recently published in the *Mathematical Spectrum* journal (see reference [5]).

    In section 4 we also establish Lemma 2 (lowest terms lemma). Using Lemma 2, in conjunction with the parametric formulas of Section 6, enables us to arrive at our objective in Section 7 which is a parametric description of all integral triangles with a 120° angle, and whose angle bisectors (of their 120° angles) have integer lengths.

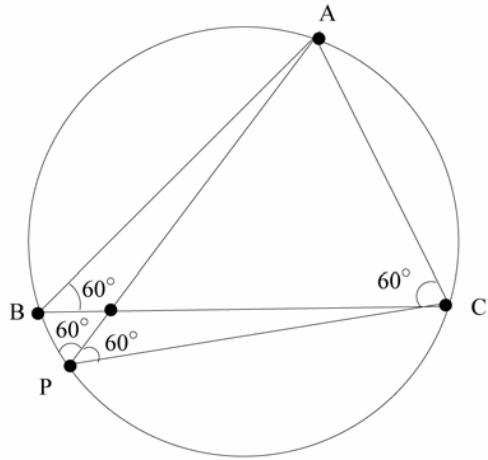

Figure 1

## 2  Ptolemy's theorem for cyclic quadrilaterals

Let *ABCD* be a cyclic quadrilateral with $\overline{AC}$ and $\overline{BD}$ being the two diagonals. Then

$$(AC) \cdot (BD) = (AB) \cdot (CD) + (CB) \cdot (AD)$$

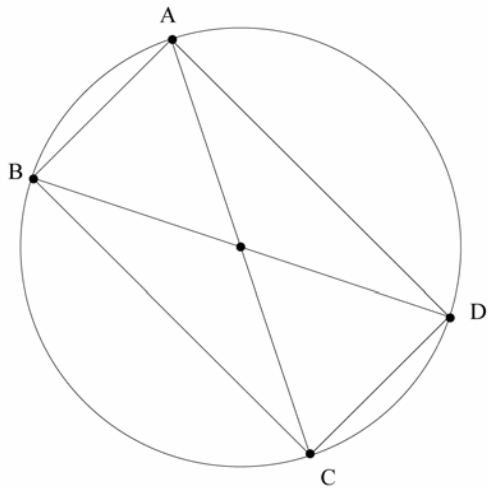

Figure 2

# 3 Proposition 1 and its proof

**Propostition 1.** *Let ABC be a triangle with the degree measure of the angle at A being 120° angle, and suppose $\overline{AD}$ is its angle bisector. Then,*

$$\frac{1}{AD} = \frac{1}{AB} + \frac{1}{AC}$$

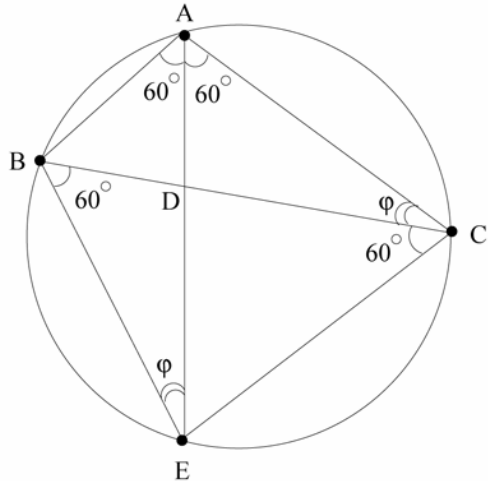

Figure 3

*Proof.* Consider the circumscribed circle of triangle *ABC* in Figure 3, and extend the angle bisector $\overline{AD}$ up to its second point of intersection with the circle at *E*. By inspection we see that ($m(\angle XYZ)$ stands for the degree measure of angle $\angle XYZ$)

$m(\angle BAE) = 60° = m(\angle DAC)$ and
$m(\angle ACB) = \varphi = m(\angle AEB)$.
Therefore, the triangles *ABE* and *ADC* are similar. We have the similarity ratios (two of the three),

$$\frac{AB}{AD} = \frac{AE}{AC} \tag{1}$$

Also, since the triangle *BCE* is equilateral,

$$BE = EC = BC \tag{2}$$

We apply Ptolemy's Theorem on the cyclic quadrilateral *ABCE*:

$$(AE) \cdot (BC) = (AB) \cdot (CE) + (AC) \cdot (BE) \qquad (3)$$

Combining (2) and (3) yields

$$AE = AB + AC \qquad (4)$$

Since $AE = AD + DE$, equality (4) implies

$$AD + DE = AB + AC,$$

and by multiplying across with $\dfrac{1}{AC}$ we arrive at

$$\frac{AD}{AC} + \frac{DE}{AC} = \frac{AB}{AC} + 1 \qquad (5)$$

From (1) and AE = AD+DE we also obtain

$$\frac{AB}{AD} = \frac{AD}{AC} + \frac{DE}{AC} \qquad (6)$$

Combining (5) with (6) implies

$$\frac{AB}{AD} = \frac{AB}{AC} + 1$$

which, in turn gives

$$\frac{1}{AD} = \frac{1}{AC} + \frac{1}{AB}$$

□

## 4    Two lemmas from number theory

The following lemma is Euclid's Lemma, a key lemma in establishing the fundamental theorem of arithmetic; and it can be found in every introductory text of number theory. For example, see [4].

**Lemma 1.** *(Euclid's lemma) if a,b,c are positive integers and a is a divisor of bc, and (a,b) =1 (i.e., a and b are relatively prime), then a must be a divisor of c.*

The lemma below will be used in Section 7.

**Lemma 2.** *(lowest terms lemma) Let a,b,c,d be positive integers such that (a,b)=1 and $\frac{a}{b} = \frac{c}{d}$. Then $c = D \cdot a$ and $d = D \cdot b$ for some positive integer D (note that D is none other than the greatest common divisor of c and d). In particular, if (c,d)=1 then a=c and b=d.*

*Proof.* We have *ad=bc*. Since *a* is a divisor of *bc* and (a,b) =1, it follows from Lemma 1 that *a* must be a divisor of *c*. Hence, $c = a \cdot D$, for some positive integer D. We have $ad = b \cdot a \bullet D; d = b \cdot D$. When $(c,d) = 1$, obviously $D = 1$.

# 5 The Diophantine equation $\frac{1}{z} = \frac{1}{x} + \frac{1}{y}$ and its solution set

**Propostion 2.** *All the positive integer solutions to the Diophantine equation*

$$\frac{1}{z} = \frac{1}{x} + \frac{1}{y} \tag{7a}$$

*are given by the parametric formulas*

$$x = k(m+n)m,\ y = k(m+n)n,\ z = kmn \text{ where } k,m,n \tag{7b}$$

*are any positive integers such that (m,n)=1*

*Proof.* If *x=km(m+n), y=kn(m+n), z=kmn* for positive integers *m,n,k* then a straightforward calculation shows that the triple {*x,y,z*} is the solution to equation (7a).

Now, let us prove the converse. If {*x,y,z*} is a positive integer solution to (7a), then *x,y,z* must have the required form given by (7b). Indeed, let *d=(x,y)*, the greatest common divisor of *x* and *y*. We have

$$\begin{cases} x = d \cdot m,\ y = d \cdot n \text{ for} \\ \text{some positive integers } m \text{ and } n \text{ with } (m,n) = 1 \end{cases} \tag{8}$$

From (7a) we obtain *z(x+y)=xy* and by (8) we further get

$$z(m+n) = dmn \tag{9}$$

Since $(m,n) = 1$, it easily follows that

$$(m+n,\ mn) = 1 \tag{10}$$

Combining (9), (10) and Lemma 1, we deduce that *mn* must be a divisor of *z* and *m+n* a divisor of *d*. Thus,

$$z = kmn \tag{11}$$

for some positive integer *k*. From (11) and (9) we further obtain *d=k(m+n)* which, in turn, implies by (8) that *x = km(m+n)* and *y=kn(m+n)*. The proof is complete.

□

# 6     Integral Triangles with a 120° angle

The parametric formulas (and their derivation) stated below can be found in reference [5]. They are part of an article published in the *Mathematical Spectrum,* Volume 43, 2010/2011, Number 2.

The entire set of the 120° triples (*a,b,c*) can be thought is as the non disjoint or overlapping union of four groups or families of integral triangles with 120° angles. (Note: *a,b,c* are the side lengths with *c* being the length of the side across from the 120° triangle.)

These four groups or families are:

$$F_1 \quad \begin{cases} a = \dfrac{d(-2rt + r^2 - 3t^2)}{4}, b = drt, \ c = \dfrac{d(r^2 + 3t^2)}{4} \\ \text{with } d, r, t \in \mathbb{Z}^+, (r,t) = 1, d \equiv 0 \ (\text{mod } 4), \\ r + t \equiv 1 (\text{mod } 2) \text{ and } r > 3t \end{cases} \tag{12}$$

$$F_2 \quad \begin{cases} a = \dfrac{d(-2rt + 3t^2 - r^2)}{4}, b = drt, \ c = \dfrac{d(r^2 + 3t^2)}{4} \\ \text{with } d, r, t \in \mathbb{Z}^+, (r,t) = 1, d \equiv 0 \ (\text{mod } 4), \\ r + t \equiv 1 (\text{mod } 2) \text{ and } r < t \end{cases} \tag{13}$$

$$F_3 \quad \begin{cases} a = \dfrac{d(-2rt+r^2-3t^2)}{4}, b = drt, \ c = \dfrac{(r^2+3t^2)}{4} \\ \text{with } d,r,t \in \mathbb{Z}^+, (r,t)=1, r \equiv t \equiv 1 \ (\text{mod } 2), \\ r > 3t \end{cases} \quad (14)$$

$$F_4 \quad \begin{cases} a = \dfrac{d(-2rt+3t^2-r^2)}{4}, b = drt, \ c = \dfrac{d(r^2+3t^2)}{4} \\ \text{with } d,r,t \in \mathbb{Z}^+, (r,t)=1, r \equiv t \equiv 1 \ (\text{mod } 2), \\ r < t \end{cases} \quad (15)$$

## 7 When the bisector of the 120° angle is of integral length

In this last section, we parametrically describe the set of all integral triangles with a 120° angle, which, also has the bisector of the 120° angle of integral length. Let $ABC$ be such a triangle with $C$ being the vertex of the angle; that is, $m(\angle ACB)=120°$ with side (integer) lengths $AC=b$, $CB=a$, $BA=c$, and with the bisector $\overline{CD}$ of the 120° angle having integer length $z$. By Proposition 1 we have,

$$\frac{1}{z} = \frac{1}{a} + \frac{1}{b} \quad (16)$$

Since $a,b,z$ are positive integers, it follows from Proposition 2 that we must also have

$$\begin{cases} a = km(m+n), b = kn(m+n), z = kmn \\ \text{for positive integers } k,m,n \text{ with } m \text{ and } n \text{ being} \\ \text{raltively prime integers; } (m,n)=1 \end{cases} \quad (17)$$

We combine (17) with the fact that the 120° integral triangle triple $(a,b,c)$ must belong to one of the families $F_1$, $F_2$, $F_3$, or $F_4$, i.e., the triple $(a,b,c)$ must satisfy one of (12), (13), (14), or (15).

Regardless of which one is satisfied, we see that in all cases we must have

$$\frac{a}{b} = \frac{m}{n} \quad (18)$$

Moreover, if the triple $(a,b,c)$ belongs to family $F_1$ or family $F_3$, then by (12) or (14) respectively, combined with (18), we get

$$\frac{m}{n} = \frac{-2rt + r^2 - 3t^2}{4rt} \qquad (19)$$

While is *(a,b,c)* belongs to family $F_2$ or family $F_4$, then by (13) or (15) respectively, combined with (18), we get

$$\frac{m}{n} = \frac{-2rt + 3t^2 - r^2}{4rt} \qquad (20)$$

Consider the two fractions (in terms of *r* and *t*) on the left-hand sides of (19) and (20). If the numerator and denominator in each fraction are not relatively prime, then they must have a prime divisor *p* in common. In view of the fact that $(r,t)=1$, one easily sees that in both cases of (19) and (20), $p=3$. This shows that in each case the greatest common divisor $(-2rt + r^2 - 3t^2, 4rt)$ is either equal to 1 or otherwise a power of 3. Likewise, for the greatest common divisor $(-2rt + 3t^2 - r^2, 4rt)$. But, neither greatest common divisor can be divisible by 9. Here is why. Take the case of (19) and suppose that both $-2rt + r^2 - 3t^2$ and $4rt$ are divisible by 9. First observe that by virtue of $(r,t)=1$, only *r* can be divisible by 3, while *t* is not divisible by 3. So then *r* would be divisible by 9. But, then $-2rt + r^2 - 3t^2$ would be exactly divisible by 3 (i.e., the highest power of 3 dividing would be 1); and so $-2rt + r^2 - 3t^2$ would not be divisible by 9. We conclude that the greatest common divisor $(-2rt + 3t^2 - r^2, 4rt)$ can only by divisible by 3 (when $r \equiv 0 \pmod 3$) but not by 9.

We have, in effect, shown that either

$$(-2rt + r^2 - 3t^2, 4rt) = 1 = (-2rt + 3t^2 - r^2, 4rt) \qquad (21)$$

or alternatively,

$$(-2rt + r^2 - 3t^2, 4rt) = 3 = (-2rt + 3t^2 - r^2, 4rt) \qquad (22)$$

Now, we apply Lemma 2 (lowest terms lemma). If (21) holds then, in conjunction with (19), (20) and Lemma 2, we get

$$\begin{cases} m = -2rt + r^2 - 3t^2, n = 4rt \\ \text{when } (a,b,c) \text{ belongs to family } F_1 \text{ or family } F_3 \end{cases} \qquad (23)$$

and

$$\begin{cases} m = -2rt + 3t^2 - r^2, n = 4rt \\ \text{when } (a,b,c) \text{ belongs to family } F_2 \text{ or } F_4 \end{cases} \qquad (24)$$

If (22) holds then, in conjunction with (19), (20), and lemma 2, we obtain

$$\begin{cases} m = 3(-2rt + r^2 - 3t^2), n = 12rt \\ \text{when } (a,b,c) \text{ belongs to family } F_1 \text{ or } F_3 \end{cases} \quad (25)$$

and

$$\begin{cases} m = 3(-2rt + 3t^2 - r^2), n = 12rt \\ \text{when } (a,b,c) \text{ belongs to family } F_2 \text{ or } F_4 \end{cases} \quad (26)$$

When $(a,b,c)$ belongs to family $F_1$ or $F_3$ and (21) holds, then by (23) or (24) respectively, combined with (12) or (14) and with (17) easily yields

$$d = 4k(-rt + r^2 - 3t^2) \text{ or}$$
$$d = 4k(-rt + 3t^2 - r^2) \text{ respectively}$$

Working similarly in the case of (22), combined with (25) or (26) respectively, and with (13) or (14), and with (17), a similar calculation produces

$$d = \frac{4k(-rt + r^2 - 3t^2)}{3} \text{ or}$$
$$d = \frac{4k(-rt + 3t^2 - r^2)}{3} \text{ respectively.}$$

## 8  Summary of results

A. Suppose that $(a,b,c)$ is a 120° integral triangle triple that belongs to family $F_1$ or $F_3$ as described by the parametric formulas in (12) or (14), then the bisector of the 120° angle has integral length $z$ if, and only if, the integer $d$ (in the formulas (12)) is of the form $d = 4k(-rt + r^2 - 3t^2)$ with $r$ not a multiple of 3; or, alternatively, of the form $d = \frac{4k(-rt + r^2 - 3t^2)}{3}$ with $r \equiv 0 \pmod{3}$, and $k$ (in both cases) an arbitrary positive integer. In the first case (when $r$ is not divisible by 3),

$$z = 4krt(-2rt + r^2 - 3t^2)$$

In the second case (when $r \equiv 0 \pmod 3$), $z = 36krt(-2rt + r^2 - 3t^2)$.

B. Supposed that $(a,b,c)$ is a 120° integral triangle triple that belongs to family $F_2$ or $F_4$ described in (13) or (15). Then the bisector of the 120° angle has integral length $z$ if, and only if, either

$$d = 4k(-rt + 3t^2 - r^2), r \text{ not divisible by 3; in which case}$$

$$z = 4krt(-2rt + 3t^2 - r^2)$$

or otherwise,

$$d = \frac{4k(-rt + 3t^2 - r^2)}{3}, r \equiv 0 \pmod{3} \text{ ; and in which case}$$

$$z = 36rt(-2rt + 3t^2 - r^2)$$

## References

[1] http://en.wikipedia.org/wiki/Ptolemy's Theorem.

[2] Coxeter, H.S.M and Greitzer, S.L., *Ptolemy's Theorem and It's Extensions,* Section 2.6 in Geometry Revisited, Washington, D.C., *Math Assoc. America,* pp. 42-43, (1967).

[3] www.cut-the_know.org/proofs/ptolemy.shtml.

[4] Kenneth H. Rose, *Elementary Number Theory and Its Applications,* fifth edition, Pearson, Addison Wesley, (2005)

   For Lemma 1 (Lemma 3.4 in the above book), see page 109.

[5] Konstantine Zelator, *Integral triangles with a 120• angle* Mathematical Spectrum, Volume 43, 2010/2011, Number 2, ISSN 0025-5653, Applied Probability Trust 2011, pp. 60-64.